\documentclass[11pt]{article}
\usepackage[left=1in,top=1in,right=1in,bottom=1in,head=.1in,nofoot]{geometry}

\usepackage{setspace,url,bm,amsmath, listings} 

\usepackage{bbm}
\usepackage{latexsym}
\usepackage{caption}
\usepackage[margin=20pt]{subcaption}
\usepackage{hyperref}
\usepackage{enumerate}

\newcommand{\GG}[1]{}

\usepackage{amsthm}
\usepackage{comment}
\theoremstyle{definition}

\newtheorem*{theorem*}{Theorem}
\newtheorem{theorem}{Theorem}

\newtheorem{lemma}{Lemma}

\newtheorem*{corollary*}{Corollary}

\usepackage{natbib} 
\bibpunct{(}{)}{;}{a}{}{,} 

\usepackage{etoolbox} 
\apptocmd{\sloppy}{\hbadness 10000\relax}{}{} 

\begin{document}
\onehalfspacing 
\title{\bf \Large The Frisch--Waugh--Lovell Theorem for Standard Errors}
\author{Peng Ding
\footnote{
Department of Statistics, University of California, Berkeley. Email: \url{pengdingpku@berkeley.edu}. Address: 425 Evans Hall, Berkeley CA 94720, USA. 
}
}
\date{}
\maketitle

\begin{abstract}
The Frisch--Waugh--Lovell Theorem states the equivalence of the coefficients from the full and partial regressions. I further show the equivalence between various standard errors. Applying the new result to stratified experiments reveals the discrepancy between model-based and design-based standard errors. 

\medskip 
\noindent {\bf Keywords:} autocorrelation; covariance estimator; clustering; heteroskedasticity; partial regression; stratified experiments
\end{abstract}

\section{Introduction}

The Frisch--Waugh--Lovell Theorem is a standard result in most econometrics textbooks, although it is less common to appear formally in statistics textbooks on linear models. The theorem states that the coefficient of $X_2$ in the full ordinary least squares (OLS) fit of $Y$ on $X = (X_1, X_2)$ equals the coefficient of $\tilde{X}_2$ in the partial OLS fit of $\tilde{Y}$ on $\tilde{X}_2$, where $\tilde{Y}$ and $\tilde{X}_2$ are the residuals from the OLS fits of $Y$ and $X_2$ on $X_1$, respectively. It is rare to see further discussions on the standard errors of the coefficients based on these two OLS fits. I shall show that the standard errors assuming homoskedasticity differ only due to a correction for degrees of freedom, and more interestingly, the Eicker--Huber--White (EHW) standard errors \citep{mackinnon1985some, angrist2008mostly} allowing for heteroskedasticity are the same in full and partial OLS fits.

I shall further extend the results to the standard errors allowing for heteroskedasticity and clustering \citep{liang1986longitudinal, angrist2008mostly}, and the standard errors allowing for heteroskedasticity and autocorrelation \citep{newey1987simple, andrews1991heteroskedasticity, andrews1992improved,  newey1994automatic, lumley1999weighted}.

Finally, I shall illustrate the results with an example from the regression analysis of stratified experiments. The example will show  the discrepancy between the standard errors based on OLS and the correct standard error based on the design of experiments.


\section{OLS and EHW standard errors}\label{sec::ols-ehw}

\subsection{Main results}
Consider the linear regression $Y = X_1 \beta_1 + X_2 \beta_2 + \varepsilon$, where $Y$ is an $n\times 1$ response vector, $X_1$ is an $n\times K$ matrix,  $X_2$ is an $n\times L$ matrix, and $ \varepsilon$ is an $n\times 1$ error vector. Let $\hat{\beta}_1$ and $ \hat{\beta}_2$ be the OLS coefficients, and $\hat{ \varepsilon } = Y -  X_1 \hat{\beta}_1 - X_2 \hat{\beta}_2$ be the residual vector. The second component $ \hat{\beta}_2$ has estimated covariance $\hat{V}$ equaling the $(2,2)$th block of
$
  \hat{\sigma}^2 (X'X)^{-1} 
$ assuming homoskedasticity, or $\hat{V}_{\textsc{ehw}}$ equaling the $(2,2)$th block of
$
  (X'X)^{-1} X' \hat{\Omega}_\textsc{ehw} X (X'X)^{-1}
$ 
allowing for heteroskedasticity, 
where $\hat{\sigma}^2 = \| \hat{ \varepsilon} \|_2^2/(n-K-L)$ is the unbiased estimator of the common variance of the errors and $\hat{\Omega}_\textsc{ehw} = \text{diag}\{ \hat{ \varepsilon}^2\}$ is the diagonal matrix consisting of the squared residuals. 
 
Consider the following partial regression. First, from the OLS fit of $Y$ on $X_1$ we obtain the residual vector $\tilde{Y}$. Then, from the column-wise OLS fit of $X_2$ on $X_1$ we obtain the residual matrix $\tilde{X}_2$. Finally, from the OLS fit of $\tilde{Y}$ on $\tilde{X}_2$ we obtain the coefficient $\tilde{\beta}_2$ and the residual vector $\tilde{\varepsilon} = \tilde{Y} -\tilde{X}_2  \tilde{\beta}_2 $. Then $\tilde{\beta}_2$ has estimated covariance 
$
\tilde{V} =   \tilde{\sigma}^2 ( \tilde{X}_2 ' \tilde{X}_2)^{-1} 
$ 
assuming homoskedasticity, or 
$
\tilde{V}_{\textsc{ehw}} = ( \tilde{X}_2 ' \tilde{X}_2)^{-1}  \tilde{X}_2'  \tilde{\Omega}_\textsc{ehw} \tilde{X}_2  ( \tilde{X}_2 ' \tilde{X}_2)^{-1} 
$
 allowing for heteroskedasticity, where $\tilde{\sigma}^2 = \| \tilde{ \varepsilon} \|_2^2/(n-L)$ and $\tilde{\Omega}_\textsc{ehw} = \text{diag}\{ \tilde{ \varepsilon}^2\}$. 

The following theorem is well known.

\begin{theorem}\label{thm::fwl}
$ \hat{\beta}_2 = \tilde{\beta}_2$. 
\end{theorem}

Econometricians often formally call it the Frisch--Waugh--Lovell (FWL) Theorem although it is also known in statistics ever since \citet{yule1907theory}. See \citet[][Section 2.5.6]{agresti2015foundations} for a textbook discussion. 

Few textbooks comment on the estimated covariances based on the full and partial regressions. In fact, we have the following theorem.

\begin{theorem}\label{thm::fwl-se}
$(n-K-L)\hat{V} = (n-L)  \tilde{V} $ and $\hat{V}_{\textsc{ehw}} = \tilde{V}_{\textsc{ehw}}.$
\end{theorem}

\citet[][Section 1]{amemiya1985advanced} and \citet[][Chapter 3]{greene2011econometric} gave elegant proofs of Theorem \ref{thm::fwl}. \citet[][page 1002--1003]{lovell1963seasonal} highlighted the importance of correcting for degrees of freedom in the partial regression, and \citet[][Chapter 1, Problem 4]{hayashi2000econometrics} hinted at the first result in Theorem \ref{thm::fwl-se}. I am not aware of a formal proof of the second result in Theorem \ref{thm::fwl-se}, and I shall prove both theorems for completeness. Two  key facts ensure $\hat{V}_{\textsc{ehw}} = \tilde{V}_{\textsc{ehw}}$: first, the residual vectors from the full and partial regressions are the same; second, the $(2,1)$th block of $(X'X)^{-1} X$ equals $( \tilde{X}_2'\tilde{X}_2 )^{-1}\tilde{X}_2'$. See the Appendix for more details.

It is important to note that both Theorems \ref{thm::fwl} and \ref{thm::fwl-se} are pure linear algebra facts without any probabilistic assumptions on the data generating process.

\subsection{Remarks}

Another known form of the FWL Theorem is that $\hat{\beta}_2  = ( \tilde{X}_2'\tilde{X}_2 )^{-1}  \tilde{X}_2'   Y$, which also follows from \eqref{eq::fwl-proof} in the Proof of Theorem \ref{thm::fwl}. So $\hat{\beta}_2$ equals the coefficient of $\tilde{X}_2$ in the OLS fit of $Y$ on $\tilde{X}_2$. While this OLS fit can give the correct point estimator, the corresponding variance estimators are not the same as the correct one even after adjusting for degrees of freedom. This is because the residual vector from this regression does not equal the correct residual vector from the full regression.

The second result of Theorem \ref{thm::fwl-se} holds only for the original form of the EHW (also called ``HC0''), but not for other variants, for example, ``HC1,'' ``HC2,'', ``HC3,'' and ``HC4'' \citep{mackinnon1985some, angrist2008mostly, zeileis2004econometric}. HC1 has a correction for degrees of freedom, so the corresponding variance estimates are $\hat{V}_{\textsc{ehw}}\times n/(n-K-L)$ and $\tilde{V}_{\textsc{ehw}}\times n/(n-L)$, which satisfy the same relationship as the covariance estimators assuming homoskedasticity. However, other variants depend on the leverage scores of the OLS fits, and there is no simple relationship between the HC$j\ (j=2,3,4)$ obtained from the full and partial regressions. The equivalence of the EHW covariance estimators (HC0) obtain by these two regressions can also be a problem especially when the dimension of $X_1$ is high and the correction for degrees of freedom becomes important \citep[e.g.,][]{baltagi2008econometric}.

\section{Extension to other robust standard errors} \label{sec::otherrobustse}
 
 From $\hat{V}_\textsc{ehw} = \tilde{V}_\textsc{ehw}$ and its proof,   a general equivalence relationship holds as long as the $\hat{\Omega}$ matrix depends only on the residual vector. However, it does not hold if $\hat{\Omega}_\textsc{ehw}$ (or $\tilde{\Omega}_\textsc{ehw}$)  depends on covariate matrix $X$ (or $\tilde{X}_2$) including its dimension, leverage scores, etc. This observation allows us to extend Theorem \ref{thm::fwl-se} to other robust standard errors. I give two important examples.

\subsection{Liang--Zeger cluster-robust standard error} \label{sec::lz}
 
When observations are clustered, for example, by states, villages, or classrooms, we can still use the OLS estimator, but it is more reasonable to use the cluster-robust variance estimator. Assume the observations are sorted based on the cluster indicators, and let $\hat{\varepsilon}_g$ be the residual vector for the observations in cluster $g\ (g=1,\ldots, G)$. Based on the full regression, $\hat{\beta}_2$ has cluster-robust covariance estimator $\hat{V}_\textsc{lz}$ \citep{liang1986longitudinal, angrist2008mostly, cameron2015practitioner} equaling the $(2,2)$th block of
$
(X'X)^{-1} X'   \hat{\Omega}_\textsc{lz}  X (X'X)^{-1} ,
$
where $ \hat{\Omega}_\textsc{lz} =  \text{diag}\{  \hat{\varepsilon}_1 \hat{\varepsilon}_1', \ldots, \hat{\varepsilon}_G \hat{\varepsilon}_G'    \} $ is the block diagonal matrix of the outer products of the residual vectors within clusters. The partial regression gives point estimator $\tilde{\beta}_2$, residual vectors $\tilde{\varepsilon}_g \ (g=1,\ldots, G)$, and cluster-robust covariance estimator equaling $\tilde{V}_\textsc{lz}
=
( \tilde{X}_2 ' \tilde{X}_2)^{-1} \tilde{X}_2'  \tilde{\Omega}_\textsc{lz} \tilde{X}_2 (\tilde{X}_2' \tilde{X}_2)^{-1} ,
$
 where 
 $\tilde{\Omega}_\textsc{lz} = \text{diag}\{  \tilde{\varepsilon}_1 \tilde{\varepsilon}_1', \ldots, \tilde{\varepsilon}_G \tilde{\varepsilon}_G'    \}   $.

We can easily extend Theorem \ref{thm::fwl-se} and show that these two cluster-robust covariance estimators are identical numerically.

\begin{theorem}
\label{thm::fwl-liangzeger}
$\hat{V}_\textsc{lz} = \tilde{V}_\textsc{lz}$.
\end{theorem}

\citet{cameron2015practitioner} stated but did not give a proof for a special case of Theorem \ref{thm::fwl-liangzeger} in panel data regression, where $X_1$ is the matrix consisting of the fixed effects dummy variables. In fact, Theorem \ref{thm::fwl-se} works for general full and partial regressions, and I shall give a formal proof in the Appendix.

Similar to the EHW covariance estimator, the LZ covariance estimator has several variants \citep{cameron2015practitioner}. Sometimes, it is further multiplied by $G/(G-1)$. This will not change Theorem \ref{thm::fwl-liangzeger}. Sometimes, it is further multiplied by $G/(G-1)\times (n-1)/(n-K-L)$ in the full regression and by $G/(G-1)\times (n-1)/(n-L)$ in the partial regression. In this case, Theorem \ref{thm::fwl-liangzeger} holds up to a correction for degrees of freedom. Other variants depend on the leverage scores \citep{mccaffrey2002bias, angrist2008mostly}, and do not satisfy Theorem \ref{thm::fwl-liangzeger} in general.

\subsection{Heteroskedasticity and autocorrelation robust standard error}\label{sec::hac}

The consistency and asymptotic Normality of the OLS estimator hold even when the error terms are serially correlated \citep{hansen1982large}. 
However, we must correctly estimate the asymptotic covariance matrix for statistical inference. 
Assume that the observations  are ordered so that $i=1,\ldots, n$ index the time. Based on the full regression, $\hat{\beta}_2$ has heteroskedasticity and autocorrelation (HAC) robust covariance estimator $\hat{V}_\textsc{hac}$   equaling the $(2,2)$th block of 
$
(X'X)^{-1} X'   \hat{\Omega}_\textsc{hac}  X (X'X)^{-1} ,
$
where $\hat{\Omega}_\textsc{hac} = (  \hat{w}_{|i-j|}   \hat{\varepsilon}_i \hat{\varepsilon}_j  )_{1\leq i,j \leq n}$ have many different forms depending on the choice of the weight $\hat{w}_{|i-j|}$. Based on the partial regression, the HAC robust covariance estimator $\tilde{V}_\textsc{hac}$   equals $( \tilde{X}_2 ' \tilde{X}_2)^{-1} \tilde{X}_2'  \tilde{\Omega}_\textsc{hac} \tilde{X}_2 (\tilde{X}_2' \tilde{X}_2)^{-1}$, where $\tilde{\Omega}_\textsc{hac} = (  \tilde{w}_{|i-j|}   \tilde{\varepsilon}_i \tilde{\varepsilon}_j  )_{1\leq i,j \leq n}$.

 \begin{theorem}
\label{thm::fwl-hac}
If $ \hat{w}_{|i-j|}  = \tilde{w}_{|i-j|}$, then 
$\hat{V}_\textsc{hac} = \tilde{V}_\textsc{hac}$. The condition holds if $\hat{w}_{|i-j|}  $ (or $ \tilde{w}_{|i-j|}$) depends on the residual vector only but not  $X$ (or $\tilde{X}_2$). 
\end{theorem}

The condition in Theorem \ref{thm::fwl-hac} holds for \citet{newey1987simple}'s estimator with a fixed lag, and \citet{lumley1999weighted}'s weighted empirical adaptive variance estimator. It is common to adjust for the degrees of freedom, so Theorem \ref{thm::fwl-hac} holds up to a minor modification. We have $(n-K-L)\hat{V}_\textsc{hac} = (n-L) \tilde{V}_\textsc{hac}$ for those estimators, for example, the default choice of the \texttt{newey} command in \texttt{stata}.  
When the weight depends on the features of the covariate matrix \citep[e.g.,][]{andrews1991heteroskedasticity, andrews1992improved, newey1994automatic}, Theorem \ref{thm::fwl-hac} does not hold even after correcting for the degrees of freedom. \citet{zeileis2004econometric}'s \texttt{sandwich} package in \texttt{R} implements many HAC robust covariance estimators.

\section{Application to the regression analysis of stratified experiments}\label{sec::sre}

In this section, I shall apply Theorems \ref{thm::fwl} and \ref{thm::fwl-se} to the regression analysis of stratified experiments. 
Consider an experiment stratified on a discrete covariate $X_i\in \{1, \ldots, K\}$, with a binary treatment $Z_i$ and an outcome $Y_i$, for unit $i \in \{ 1,\ldots, n \}$. Let $[k]$ be the set of indices, $n_{[k]}$ be the sample size, and $\pi_{[k]} = n_{[k]}/n$ be the proportion of units within stratum $k$. Within stratum $k$, let $n_{[k]1}$ and $n_{[k]0}$ be the number of units under treatment and control, respectively, and define $e_{[k]} = n_{[k]1}/n_{[k]}$ as the propensity score, i.e., the probability of receiving the treatment. Assume that $0< e_{[k]} < 1$ for all $k$.  In general, different strata can have different probabilities of receiving the treatment and therefore the $e_{[k]}$'s can vary across strata.

A question of interest is the treatment effect of $Z$ on the $Y$. A standard estimator is the coefficient of $Z_i$ in the OLS fit of $Y_i$ on the treatment indicator $Z_i$ and $K$ stratum indicators:
\begin{eqnarray}\label{eq::ols}
Y_i \sim Z_i +  I(i \in [1]) + \cdots  I( i \in [K] ), 
\end{eqnarray}
where the intercept is omitted to avoid overparametrization \citep[e.g.,][]{duflo2007using, angrist2014opportunity, ImbensRubin15, gibbons2018broken}.

\subsection{Point estimator}

 Interestingly, we can use Theorem \ref{thm::fwl} to obtain an explicit formula for this OLS estimator. 

\begin{theorem}
\label{thm::ols-fitted}
The OLS coefficient of $Z_i$ from \eqref{eq::ols} equals 
$
\hat{\tau}_{\textsc{ols}}
=    \sum_{k=1}^K   \omega_{[k]}   \hat{\tau}_{[k]}   
$
where $\omega_{[k]} =  \pi_{[k]}  e_{[k]} (1-e_{[k]}) / \sum_{k'=1}^K  \pi_{[k']}  e_{[k']} (1-e_{[k']})$ 
and
$\hat{\tau}_{[k]}  = \bar{Y}_{[k]1} - \bar{Y}_{[k]0} $.
\end{theorem}

We can prove Theorem \ref{thm::ols-fitted} by directly solving the Normal equation of OLS in \eqref{eq::ols}. But using Theorem \ref{thm::fwl} yields a more elegant proof. I give a sketch below and present the details in the supplementary material. 

\begin{proof}
From the OLS fit of $Y_i$ on the stratum indicators, we obtain the residuals $\tilde{Y}_i  = Y_i- \bar{Y}_{[k]}$ for unit with $i \in [k]$, where 
$
 \bar{Y}_{[k]} = n_{[k]}^{-1} \sum_{i \in [k]} Y_i
 = e_{[k]}   \bar{Y}_{[k]1} + (1 - e_{[k]}) \bar{Y}_{[k]0}.
$
From the OLS fit of $Z_i$ on the stratum indicators, we obtain residuals $\tilde{Z}_i = Z_i - e_{[k]}$ for unit with $i \in [k]$. Using Theorem \ref{thm::fwl}, we have $ \hat{\tau}_{\textsc{ols}} = \sum_{i=1}^n  \tilde{Z}_i \tilde{Y}_i  /  \sum_{i=1}^n  \tilde{Z}_i ^2$ based on a one-dimensional OLS. The conclusion follows by simplifying the denominator and the numerator: the former reduces to
$n \sum_{k=1}^K \pi_{[k]}   e_{[k]}(1- e_{[k]})$,
and the latter reduces to 
$ n \sum_{k=1}^K \pi_{[k]}   e_{[k]}(1- e_{[k]}) \hat{\tau}_{[k]}.
$
 \end{proof}

Theorem \ref{thm::ols-fitted} states that $\hat{\tau}_{\textsc{ols}}$ is a weighed average of the difference-in-means of the outcomes within stratum $k$, with weights proportional to  the product of $\pi_{[k]} $ and the variance of the treatment $ e_{[k]} (1-e_{[k]})$. It is a numerical result of the OLS fit, which holds without any probabilistic assumptions on the data generating process. \citet{angrist1998estimating} and \citet[][Section 3.3]{angrist2008mostly}, \citet[][Chapter 9]{ImbensRubin15}, \citet[][Theorem A.1]{aronow2016does}, and \citet[][Proposition 1]{gibbons2018broken} gave similar forms of Theorem \ref{thm::ols-fitted}, assuming certain superpopulation sampling and asymptotic scheme.

\subsection{Variance estimators}

Within stratum $k$, we have a completely randomized experiment, and therefore we can estimate the variance of $\hat{\tau}_{[k]}$ by $\hat{V}_{[k]} = s_{[k]1}^2 / n_{[k]1} +  s_{[k]0}^2 / n_{[k]0}$, where $s_{[k]1}^2  = n_{[k]1}^{-1} \sum_{i\in [k], Z_i=1} (Y_i - \bar{Y}_{[k]1})^2$ and $s_{[k]0}^2  = n_{[k]0}^{-1} \sum_{i\in [k], Z_i=0} (Y_i - \bar{Y}_{[k]0})^2$ are the sample variances of the outcomes  under treatment and control, respectively. The usual definitions of the sample variances have denominators $ n_{[k]1}-1$ and $ n_{[k]0}-1$, but we adopt these definitions to simplify the expressions below. The difference is minor when the sample sizes $n_{[k]1}$ and $ n_{[k]0}$ are large. 
The variance estimator $\hat{V}_{[k]}$ is the standard choice assuming either a superpopulation or a finite population \citep{ImbensRubin15, ding2017bridging}. Then based on Theorem \ref{thm::ols-fitted}, a variance estimator for $\hat{\tau}_{\textsc{ols}}$ is 
$
\hat{V}_0   = \sum_{k=1}^K    \omega_{[k]} ^2   \hat{V}_{[k]}. 
$
It is curious to see whether this variance estimator is similar to the ones based on the OLS fit. Unfortunately, the answer is no. We can use Theorem \ref{thm::fwl-se} to derive the following results.

\begin{theorem}\label{thm::variances-ols}
Based on the OLS fit in \eqref{eq::ols}, the variance estimator assuming homoskedasticity is
$$
\hat{V} = \frac{n}{n-K-1}   \sum_{k=1}^K \omega_{[k]} \pi_{[k]}     \left\{   \hat{\Lambda}_{[k]}
+   ( \hat{\tau}_{[k]} -  \hat{\tau}_{\textsc{ols}} )^2 / n_{[k]}  \right\}  . 
$$
and the variance estimator allowing for heteroskedasticity is
$$
\hat{V}_{\textsc{ehw}} =  \sum_{k=1}^K  \omega_{[k]}^2   \left\{ 
\hat{V}_{[k]} + \Delta_{[k]} ( \hat{\tau}_{[k]} -  \hat{\tau}_{\textsc{ols}} )^2 / n_{[k]} 
\right\} . 
$$
where $\hat{\Lambda}_{[k]}=  s_{[k]1}^2 / n_{[k]0} +  s_{[k]0}^2 / n_{[k]1}$, and 
$
\Delta_{[k]}  = e_{[k]}^{-1} (1 - e_{[k]})^{-1} - 3 \geq 1
$
because $e_{[k]} (1 - e_{[k]}) \leq 1/4$. 
\end{theorem}

I give the sketch of the proof below and present the details in the supplementary material. 

\begin{proof}
In the Proof of Theorem \ref{thm::fwl-se}, I have shown that the residuals from the full regression are the same as the residuals from the partial regression, which, for unit $i \in [k]$, is $\hat{\varepsilon}_i  = Y_i-   \bar{Y}_{[k]1}  + (1-e_{[k]}) (\hat{\tau}_{[k]}-  \hat{\tau}_{\textsc{ols}} )$ for treated unit and $\hat{\varepsilon}_i  = Y_i- \bar{Y}_{[k]0}  - e_{[k]} ( \hat{\tau}_{[k]} - \hat{\tau}_{\textsc{ols}}  )$ for control unit.  
The variance estimators with and without assuming homoskedasticity are
$
\hat{V} =   \sum_{i=1}^n \hat{\varepsilon}_i^2  / \{    (n-K-1)  \sum_{i=1}^n  \tilde{Z}_i ^2  \} 
$
and
$
\hat{V}_{\textsc{ehw}} =  \sum_{i=1}^n  \tilde{Z}_i ^2 \hat{\varepsilon}_i^2  /  ( \sum_{i=1}^n  \tilde{Z}_i ^2 )^2  . 
$
The conclusions follows from simplifying the denominators and numerators.
\end{proof}

The variance estimator assuming homoskedasticity is quite different from $\hat{V}_0$. The label switching of $n_{[k]1}$ and $n_{[k]0}$ in $\hat{V}_{[k]} $ and $\hat{\Lambda}_{[k]}$ appeared in the regression analysis of completely randomized experiments \citep{Ding:2014}. In general, inference based on $\hat{V} $ can be either conservative or anti-conservative. The variance estimator allowing for heteroskedasticity is close to $\hat{V}_0$ when there is no treatment effect heterogeneity across strata, i.e., $\hat{\tau}_{[k]} -  \hat{\tau}_{\textsc{ols}}$ is close to zero. Otherwise, $\hat{V}_{\textsc{ehw}}$ is more conservative than $\hat{V}_0$  because the positive term $\Delta_{[k]} ( \hat{\tau}_{[k]} -  \hat{\tau}_{\textsc{ols}} )^2 / n_{[k]} $ has the same order as  $\hat{V}_{[k]}$ when $n_{[k]1}$ and $n_{[k]0}$ are large.

\subsection{Remarks}

Theorem \ref{thm::ols-fitted} shows that the coefficient of the treatment  in the OLS fit \eqref{eq::ols} is a weighted average of the within-stratum difference-in-means estimators, with weights proportional to $ \pi_{[k]}  e_{[k]} (1-e_{[k]})$ rather than $\pi_{[k]}.$ If the weights were $\pi_{[k]}$, then the estimator would be unbiased for the usual average treatment effect for the whole population. Theorem \ref{thm::ols-fitted} suggests that $\hat{\tau}_\textsc{ols}$ is not unbiased for the usual estimand of interest unless the propensity scores are constant. 
 
Sometimes  $\hat{\tau}_\textsc{ols}$ may estimate the parameter of interest, in the cases, for example, where the propensity scores $e_{[k]} $'s are constant, or we target a population weighted by  $\pi_{[k]}  e_{[k]} (1-e_{[k]})$. Even in these cases, from Theorem \ref{thm::variances-ols}, we recommend neither of the variance estimators based on OLS because both can be quite different from $\hat{V}_0$. The one assuming homoskedasticity $\hat{V}$ can be anti-conservative, and the one allowing for heteroskedasticity $\hat{V}_\textsc{ehw}$ are overly conservative in the presence of treatment effect heterogeneity. Recently, \citet{abadie2020sampling} highlighted the difference between model-based and design-based uncertainties.

Overall, Theorems \ref{thm::ols-fitted} and \ref{thm::variances-ols} are negative results on the regression estimator for stratified experiments, which extend \citet{freedman2008regression_b}'s critique on the regression estimator for completely randomized experiments. The first critique is on the point estimator and the second critique is on the standard errors based on OLS. \citet{miratrix2013adjusting} and \citet{gibbons2018broken} discussed alternative estimators in stratified experiments, suggesting including all the interaction terms of the treatment and stratum indicators in the OLS fit \eqref{eq::ols}. This is also a strategy proposed by \citet{lin2013agnostic} in analyzing completely randomized experiments with covariates. With this fully interacted regression, the EHW standard error is nearly identical to the correct standard error.

Finally, because Theorems \ref{thm::ols-fitted} and \ref{thm::variances-ols} are purely numeric results without assuming any probabilistic model, they also hold for other studies including completely randomized experiment and observational studies with a discrete covariate $X$. \citet{lin2013agnostic} and  \citet{miratrix2013adjusting} discussed the former; \citet{angrist2008mostly} discussed the latter.

\section{Discussion}\label{sec::discuss}

The Frisch--Waugh--Lovell Theorem is a powerful tool to understand regression coefficients from full and partial regressions. This note provides some further results on the associated covariance estimators that assume homoskedasticity and allow for heteroskedasticity, clustering, and autocorrelation. 

I illustrate these results with an example from the regression analysis of stratified experiments, showing that the commonly-used OLS estimator may not be consistent for the targeted parameter of interest, and moreover, even if it is consistent under strong assumptions, the corresponding variance estimators can be either conservative or anti-conservative, assuming homoskedasticity or not. These results are also useful for understanding many other regression estimators in randomized experiments. I will report additional results in a longer paper.

\appendix 
\section*{Appendix: Proofs of the results on the robust standard errors}

Define $H = X(X'X)^{-1} X'$ and $H_1 = X_1 (X_1'X_1)^{-1} X_1 '$ as the hat matrices generated by $X$ and $X_1$, respectively. Then $\tilde{Y} = (I-H_1)Y$ and $\tilde{X}_2 =(I-H_1) X_2$ are the residuals after regressing on $X_1$. Define $\tilde{H}_2 = \tilde{X}_2 (\tilde{X}_2'\tilde{X}_2)^{-1} \tilde{X}_2 '$ as the hat matrix generated by $\tilde{X}_2 .$

The first lemma follows from the inverse of a block matrix, and I omit the proof. 

\begin{lemma}\label{lemma::block-inv}
We have
$$
(X'X)^{-1} = 
\begin{pmatrix}
X_1'X_1 & X_1'X_2\\
X_2'X_1 & X_2'X_2
\end{pmatrix}^{-1}
=
\begin{pmatrix}
S_{11} & S_{12} \\
S_{21} & S_{22}
\end{pmatrix} ,
$$
where $
S_{11} = (X_1'X_1)^{-1} + (X_1'X_1)^{-1}  X_1' X_2  ( \tilde{X}_2'\tilde{X}_2 )^{-1} X_2' X_1 (X_1'X_1)^{-1} ,
$
$
S_{21} = S_{12}' = -  ( \tilde{X}_2'\tilde{X}_2 )^{-1} X_2'X_1 (X_1'X_1)^{-1}  ,
$
and
$
S_{22} = ( \tilde{X}_2'\tilde{X}_2 )^{-1}.
$
\end{lemma}

The second lemma decomposes $H$ into two orthogonal projection matrices. It
follows from Lemma \ref{lemma::block-inv} and the definitions of the hat matrices, and I omit the proof.

\begin{lemma}
\label{lemma::hat}
$H = H_1 + \tilde{H}_2$ and $H_1 \tilde{H}_2 = 0$. 
\end{lemma}

\paragraph{Proof of Theorem \ref{thm::fwl}}
From the definition of the OLS coefficients, we have
$$
\begin{pmatrix}
\hat{\beta}_1\\
\hat{\beta}_2
\end{pmatrix}
=
(X'X)^{-1}X'Y
=
\begin{pmatrix}
S_{11} & S_{12} \\
S_{21} & S_{22}
\end{pmatrix}
\begin{pmatrix}
X_1'Y \\
X_2'Y
\end{pmatrix}.
$$
Therefore, the second component is
\begin{eqnarray}
\label{eq::fwl-proof}
\hat{\beta}_2 = S_{21} X_1'Y  +S_{22} X_2'Y
= -  ( \tilde{X}_2'\tilde{X}_2 )^{-1} X_2' H_1 Y + ( \tilde{X}_2'\tilde{X}_2 )^{-1}X_2'Y
=( \tilde{X}_2'\tilde{X}_2 )^{-1}  X_2'  (I-H_1) Y.
\end{eqnarray}
Because $I-H_1$ is a projection matrix, we have $X_2'  (I-H_1) Y = X_2'  (I-H_1)^2 Y = \tilde{X}_2' \tilde{Y}$. Thus, $\hat{\beta}_2 = \tilde{\beta}_2.$

\paragraph{Proof of Theorem \ref{thm::fwl-se}}\label{sec::varianceproof}

Let $\hat{ \varepsilon} = (I-H)Y$ and $\tilde{ \varepsilon} = (I-\tilde{H}_2) \tilde{Y}$ be the residual vectors from the full and partial regressions. I first show that $\hat{ \varepsilon}  = \tilde{ \varepsilon}$. From $\tilde{ \varepsilon} = (I-\tilde{H}_2) \tilde{Y} = (I- \tilde{H}_2)(I-H_1) Y $, it suffices to show that $I-H = (I- \tilde{H}_2)(I-H_1)$, or, equivalently, $ I-H =  I - H_1 - \tilde{H}_2 + \tilde{H}_2 H_1$. This holds because of Lemma \ref{lemma::hat}. 
Therefore, we also have $\hat{\Omega}_\textsc{ehw} = \tilde{\Omega}_\textsc{ehw}$, where $ \hat{\Omega}_\textsc{ehw} = \textup{diag}\{ \hat{ \varepsilon}^2  \}$ and $\tilde{\Omega}_\textsc{ehw} =  \textup{diag}\{ \tilde{ \varepsilon}^2  \}.$

Let $\hat{\sigma}^2 = \| \hat{ \varepsilon} \|^2/(n-K-L)$ and $\tilde{\sigma}^2 = \| \tilde{ \varepsilon} \|^2/(n-L)$ be the common variance estimators. Under homoskedasticity, the covariance estimator for $\hat{\beta}_2$ is the $(2,2)$th block of $\hat{\sigma}^2 (X'X)^{-1}$, that is, $ \hat{\sigma}^2   S_{22} = \hat{\sigma}^2  ( \tilde{X}_2'\tilde{X}_2 )^{-1}$ based on Lemma \ref{lemma::block-inv}, which is identical to the covariance estimator for $\tilde{\beta}_2$ up to a correction for degrees of freedom.

The EHW covariance estimator from the full regression is the $(2,2)$ block of
$
 A\hat{\Omega}_\textsc{ehw} A',
$
where  
$$
A = (X'X)^{-1} X' = \begin{pmatrix}
*\\
 -  ( \tilde{X}_2'\tilde{X}_2 )^{-1} X_2' H_1 + ( \tilde{X}_2'\tilde{X}_2 )^{-1}X_2'
\end{pmatrix}
=\begin{pmatrix}
*\\
 ( \tilde{X}_2'\tilde{X}_2 )^{-1}\tilde{X}_2'
\end{pmatrix}.
$$
The $*$ term does not affect the final calculation.
Define $\tilde{A}_2 =   ( \tilde{X}_2'\tilde{X}_2 )^{-1}\tilde{X}_2'$, and then
$
\hat{V}_{\textsc{ehw}} 
= \tilde{A}_2 \hat{\Omega}_\textsc{ehw} \tilde{A}_2' 
=  \tilde{A}_2 \tilde{\Omega}_\textsc{ehw}  \tilde{A}_2',
$
which equals the EHW covariance estimator $\tilde{V}_{\textsc{ehw}}$ from the partial regression.

\paragraph{Proof of Theorem \ref{thm::fwl-liangzeger}}

First, $\hat{\Omega}_\textsc{lz} =  \tilde{\Omega}_\textsc{lz} $ because the residual vectors are the same. Then we can show that $\hat{V}_\textsc{lz} = \tilde{V}_\textsc{lz}$ following the same steps as the proof of Theorem \ref{thm::fwl-se} with minor modifications by replacing all the subscript ``EHW'' by ``LZ''.

\paragraph{Proof of Theorem \ref{thm::fwl-hac}}

First, $\hat{\Omega}_\textsc{hac} =  \tilde{\Omega}_\textsc{hac} $ because the residual vectors are the same and the weights are also the same. Then we can show that $\hat{V}_\textsc{lz} = \tilde{V}_\textsc{lz}$ following the same steps as the proof of Theorem \ref{thm::fwl-se} with minor modifications by replacing all the subscript ``EHW'' by ``HAC''.

\section*{Acknowledgments}
I thank Professor Alan Agresti for the reference of \citet{yule1907theory}, and two reviewers, Fangzhou Su, Chaoran Yu, Anqi Zhao, and Avi Feller for helpful comments. This work was supported by the U.S. National Science Foundation (grants 1713152 and 1945136).

\bibliographystyle{apalike}
\bibliography{causal}

\newpage
\setcounter{page}{1}
\pagenumbering{arabic}
\renewcommand*{\thepage}{S\arabic{page}}

\begin{center}
\bf \Large 
Online supplementary material for \\
``The Frisch--Waugh--Lovell Theorem for Standard Errors''

by Peng Ding
\end{center}
 
\bigskip 
\bigskip 

\section{Full proofs of the results on stratified experiments}

\paragraph{Proof of Theorem \ref{thm::ols-fitted}}
From the OLS fit of $Y_i$ on the stratum indicators, we obtain the residuals $\tilde{Y}_i  = Y_i- \bar{Y}_{[k]}$ for unit $i \in [k]$, where 
$$
 \bar{Y}_{[k]} = n_{[k]}^{-1} \sum_{i \in [k]} Y_i
 = e_{[k]}   \bar{Y}_{[k]1} + (1 - e_{[k]}) \bar{Y}_{[k]0}.
$$
From the OLS fit of $Z_i$ on the stratum indicators, we obtain residuals $\tilde{Z}_i = Z_i - e_{[k]}$ for unit with $i \in [k]$. Using Theorem \ref{thm::fwl}, we have 
$$
 \hat{\tau}_{\textsc{ols}} = \frac{  \sum_{i=1}^n  \tilde{Z}_i \tilde{Y}_i  }{   \sum_{i=1}^n  \tilde{Z}_i ^2} 
 $$ 
 based on a one-dimensional OLS. The conclusion follows by simplifying the denominator and the numerator: the former reduces to
\begin{eqnarray} 
\sum_{i=1}^n  \tilde{Z}_i ^2
&=& \sum_{k=1}^K \sum_{i \in [k]} (Z_i - e_{[k]})^2
=  \sum_{k=1}^K \left\{  n_{[k]1}(1- e_{[k]})^2 +  n_{[k]0}  e_{[k]}^2 \right\}  \nonumber \\
&=& \sum_{k=1}^K n_{[k]}  \left\{  e_{[k]}(1- e_{[k]})^2 +  (1- e_{[k]})  e_{[k]}^2 \right\}
=  n \sum_{k=1}^K \pi_{[k]}   e_{[k]}(1- e_{[k]}), \label{eq::denominator}
\end{eqnarray} 
and the latter reduces to 
\begin{eqnarray*}
\sum_{i=1}^n  \tilde{Z}_i \tilde{Y}_i 
&=& \sum_{k=1}^K \sum_{i \in [k]} (Z_i - e_{[k]}) Y_i 
= \sum_{k=1}^K \left\{ (1 - e_{[k]}) n_{[k]1}\bar{Y}_{[k]1}
-
 e_{[k]}  n_{[k]0} \bar{Y}_{[k]0}
 \right\}\\
 &=&   n \sum_{k=1}^K \pi_{[k]}   e_{[k]}(1- e_{[k]}) \hat{\tau}_{[k]}.
\end{eqnarray*}

\paragraph{Proof of Theorem \ref{thm::variances-ols}}
In the Proof of Theorem \ref{thm::fwl-se}, I have shown that the residuals from the full regression are the same as the residuals from the partial regression, which, for unit $i \in [k]$, is
\begin{eqnarray*}
\hat{\varepsilon}_i 
&=& Y_i- \bar{Y}_{[k]} - \hat{\tau}_{\textsc{ols}} ( Z_i - e_{[k]})\\
&=& Y_i-  e_{[k]} \bar{Y}_{[k]1}  - (1-e_{[k]}) \bar{Y}_{[k]0} -   \hat{\tau}_{\textsc{ols}}( Z_i - e_{[k]}) \\
&=&
\left\{
                \begin{array}{ll}
                Y_i-   \bar{Y}_{[k]1}  + (1-e_{[k]}) (\hat{\tau}_{[k]}-  \hat{\tau}_{\textsc{ols}} ), & (Z_i=1)\\
                Y_i- \bar{Y}_{[k]0}  - e_{[k]} ( \hat{\tau}_{[k]} - \hat{\tau}_{\textsc{ols}}  ), & (Z_i=0) . 
                \end{array}
              \right.
\end{eqnarray*}
The variance estimators are  
$$
\hat{V} =  \frac{  \sum_{i=1}^n \hat{\varepsilon}_i^2  }{   (n-K-1)  \sum_{i=1}^n  \tilde{Z}_i ^2  }
$$
 assuming homoskedasticity, and 
 $$
\hat{V}_{\textsc{ehw}} =\frac{  \sum_{i=1}^n  \tilde{Z}_i ^2 \hat{\varepsilon}_i^2  }{  \left( \sum_{i=1}^n  \tilde{Z}_i ^2 \right)^2  }
$$
allowing for heteroskedasticity. 
We can easily determine the denominators by \eqref{eq::denominator}, so we only need to simplify the numerators.  

\paragraph{(a)}
The numerator of $\hat{V}$ is $\sum_{i=1}^n \hat{\varepsilon}_i^2$, which equals
\begin{eqnarray*}
&&\sum_{k=1}^K  \left[ 
\sum_{i \in [k], Z_i=1}  \left\{ Y_i-   \bar{Y}_{[k]1}  + (1-e_{[k]}) (\hat{\tau}_{[k]}-   \hat{\tau}_{\textsc{ols}} ) \right\}^2
+
\sum_{i \in [k], Z_i=0} \left\{  Y_i- \bar{Y}_{[k]0}  - e_{[k]} ( \hat{\tau}_{[k]} -  \hat{\tau}_{\textsc{ols}} ) \right\}^2
\right] 
\\
&=& \sum_{k=1}^K  \left\{   n_{[k]1}   s_{[k]1}^2 + n_{[k]1} (1-e_{[k]})^2 (\hat{\tau}_{[k]}-   \hat{\tau}_{\textsc{ols}}  )^2
+
n_{[k]0}   s_{[k]0}^2  + n_{[k]0} e_{[k]}^2 ( \hat{\tau}_{[k]} -   \hat{\tau}_{\textsc{ols}} )^2 \right\} .
\end{eqnarray*}
Therefore, $\hat{V}\times (n-K-1)/n$ is
\begin{eqnarray*}
&&
\frac{  \sum_{k=1}^K  \left\{   n_{[k]1}   s_{[k]1}^2 + n_{[k]1} (1-e_{[k]})^2 (\hat{\tau}_{[k]}-   \hat{\tau}_{\textsc{ols}}  )^2
+
n_{[k]0}  s_{[k]0}^2  + n_{[k]0} e_{[k]}^2 ( \hat{\tau}_{[k]} -  \hat{\tau}_{\textsc{ols}}  )^2 \right\}  }
{  n^2 \sum_{k=1}^K \pi_{[k]}   e_{[k]}(1- e_{[k]})  } \\
& =  &
\frac{  \sum_{k=1}^K \pi_{[k]}^2 e_{[k]} (1-e_{[k]})  \left\{    \frac{ s_{[k]1}^2}{n_{[k]0}} 
+   (1-e_{[k]}) \frac{  (\hat{\tau}_{[k]}-   \hat{\tau}_{\textsc{ols}}  )^2}{n_{[k]}}
+
  \frac{ s_{[k]0}^2 }{n_{[k]1}} 
+    e_{[k]} \frac{ ( \hat{\tau}_{[k]} -  \hat{\tau}_{\textsc{ols}}  )^2}{n_{[k]}} \right\}  }
{   \sum_{k=1}^K \pi_{[k]}   e_{[k]}(1- e_{[k]})  } \\
&=& 
\frac{  \sum_{k=1}^K \pi_{[k]}^2  (1-e_{[k]} ) e_{[k]}  \left\{     \frac{ s_{[k]1}^2}{n_{[k]0}} 
+
  \frac{ s_{[k]0}^2 }{n_{[k]1}} 
+    \frac{ ( \hat{\tau}_{[k]} -  \hat{\tau}_{\textsc{ols}}  )^2}{n_{[k]}} \right\}  }
{   \sum_{k=1}^K \pi_{[k]}   e_{[k]}(1- e_{[k]})  } .
\end{eqnarray*}
The formula of $\hat{V} $ then follows.

\paragraph{(b)}
The numerator of $\hat{V}_{\textsc{ehw}} $ is $ \sum_{i=1}^n (Z_i - e_{[k]})^2 \hat{\varepsilon}_i^2 $, which equals
\begin{eqnarray*}
&&
\sum_{k=1}^K  \left[ 
(1 - e_{[k]})^2 \sum_{i \in [k], Z_i=1}   \left\{ Y_i-   \bar{Y}_{[k]1}  + (1-e_{[k]}) (\hat{\tau}_{[k]}-   \hat{\tau}_{\textsc{ols}}  ) \right\}^2 \right.   \\
&& \left.  \qquad \qquad \qquad \qquad +
e_{[k]}^2 \sum_{i \in [k], Z_i=0} \left\{  Y_i- \bar{Y}_{[k]0}  - e_{[k]} ( \hat{\tau}_{[k]} -  \hat{\tau}_{\textsc{ols}}  ) \right\}^2
\right]  \\
&=& \sum_{k=1}^K  \left\{  (1 - e_{[k]})^2 n_{[k]1}    s_{[k]1}^2 + n_{[k]1} (1-e_{[k]})^4 (\hat{\tau}_{[k]}-   \hat{\tau}_{\textsc{ols}}  )^2 \right. \\
&& \left.  \qquad \qquad \qquad \qquad  +
e_{[k]}^2 n_{[k]0}   s_{[k]0}^2  + n_{[k]0} e_{[k]}^4 ( \hat{\tau}_{[k]} -  \hat{\tau}_{\textsc{ols}}  )^2 \right\} \\
&=&  n^2
\sum_{k=1}^K \pi_{[k]}^2 \left\{  (1 - e_{[k]})^2 e_{[k]}^2  \frac{s_{[k]1}^2}{n_{[k]1}} 
+  e_{[k]} (1-e_{[k]})^4  \frac{ (\hat{\tau}_{[k]}-   \hat{\tau}_{\textsc{ols}}  )^2}{n_{[k]}}
\right. \\
&& \left.  \qquad \qquad \qquad  \qquad+
e_{[k]}^2 (1-e_{[k]})^2  \frac{s_{[k]0}^2}{  n_{[k]0} }  
+ (1-e_{[k]}) e_{[k]}^4  \frac{ ( \hat{\tau}_{[k]} -  \hat{\tau}_{\textsc{ols}}  )^2}{n_{[k]}} \right\}   \\
&=& n^2
\sum_{k=1}^K \pi_{[k]}^2 (1 - e_{[k]})^2 e_{[k]}^2  \left\{ 
\frac{s_{[k]1}^2}{n_{[k]1}}  +  \frac{s_{[k]0}^2}{  n_{[k]0} }  + \Delta_{[k]} \frac{ ( \hat{\tau}_{[k]} -  \hat{\tau}_{\textsc{ols}}  )^2}{n_{[k]}}
\right\},
\end{eqnarray*}
where
\begin{eqnarray*}
\Delta_{[k]} &=&  \frac{  (1 - e_{[k]})^2   }{ e_{[k]} }    + \frac{ e_{[k]}^2 }{ (1 - e_{[k]})}
=  \frac{  (1 - e_{[k]})^3 +   e_{[k]}^3 }{ e_{[k]} (1 - e_{[k]}) }   \\
&=& \frac{  (1 - e_{[k]})^2 +   e_{[k]}^2 - (1 - e_{[k]}) e_{[k]} }{ e_{[k]} (1 - e_{[k]}) }   \\
&=& \frac{  1  - 3(1 - e_{[k]}) e_{[k]} }{ e_{[k]} (1 - e_{[k]}) }   \\
&=& e_{[k]}^{-1} (1 - e_{[k]})^{-1} - 3.
\end{eqnarray*}
The formula of $\hat{V}_{\textsc{ehw}} $ then follows.

\section{Verify the main results in STATA}

I use STATA version 14.1 to verify the theoretical results. STATA automatically corrects for the degrees of freedom, so the standard errors from the partial regression equals the corresponding standard errors from the full regression multiplied by $\sqrt{71/72}$.

The following code loads the data and obtain the residuals from two partial regressions. 
\singlespacing
{\footnotesize  
\begin{verbatim}
. quiet use http://www.stata-press.com/data/r14/auto, clear

. quiet regress price displ

. predict price_res, residuals

. quiet regress weight displ

. predict weight_res, residuals
\end{verbatim}
}
\onehalfspacing

The following code verifies the conclusion about the standard errors assuming homoskedasticity.
\singlespacing
{\footnotesize 
\begin{verbatim}
. regress price weight displ

      Source |       SS           df       MS      Number of obs   =        74
-------------+----------------------------------   F(2, 71)        =     14.57
       Model |   184768050         2    92384025   Prob > F        =    0.0000
    Residual |   450297346        71  6342216.14   R-squared       =    0.2909
-------------+----------------------------------   Adj R-squared   =    0.2710
       Total |   635065396        73  8699525.97   Root MSE        =    2518.4

------------------------------------------------------------------------------
       price |      Coef.   Std. Err.      t    P>|t|     [95% Conf. Interval]
-------------+----------------------------------------------------------------
      weight |   1.823366   .8498204     2.15   0.035     .1288723     3.51786
displacement |   2.087054     7.1918     0.29   0.773    -12.25299     16.4271
       _cons |    247.907   1472.021     0.17   0.867     -2687.22    3183.034
------------------------------------------------------------------------------

. regress price_res weight_res

      Source |       SS           df       MS      Number of obs   =        74
-------------+----------------------------------   F(1, 72)        =      4.67
       Model |  29196745.7         1  29196745.7   Prob > F        =    0.0341
    Residual |   450297358        72  6254129.98   R-squared       =    0.0609
-------------+----------------------------------   Adj R-squared   =    0.0478
       Total |   479494104        73  6568412.38   Root MSE        =    2500.8

------------------------------------------------------------------------------
   price_res |      Coef.   Std. Err.      t    P>|t|     [95% Conf. Interval]
-------------+----------------------------------------------------------------
  weight_res |   1.823366   .8438982     2.16   0.034     .1410856    3.505646
       _cons |   2.70e-06   290.7151     0.00   1.000      -579.53      579.53
------------------------------------------------------------------------------

. display .8498204*sqrt(71/72)
.84389823
\end{verbatim}
}
\onehalfspacing

The following code verifies the conclusion about the EHW standard errors. 
\singlespacing
{\footnotesize 
\begin{verbatim}
. regress price weight displ, vce(robust)

Linear regression                               Number of obs     =         74
                                                F(2, 71)          =      14.44
                                                Prob > F          =     0.0000
                                                R-squared         =     0.2909
                                                Root MSE          =     2518.4

------------------------------------------------------------------------------
             |               Robust
       price |      Coef.   Std. Err.      t    P>|t|     [95% Conf. Interval]
-------------+----------------------------------------------------------------
      weight |   1.823366   .7808755     2.34   0.022     .2663445    3.380387
displacement |   2.087054   7.436967     0.28   0.780    -12.74184    16.91595
       _cons |    247.907   1129.602     0.22   0.827    -2004.455    2500.269
------------------------------------------------------------------------------

. regress price_res weight_res, vce(robust)

Linear regression                               Number of obs     =         74
                                                F(1, 72)          =       5.53
                                                Prob > F          =     0.0214
                                                R-squared         =     0.0609
                                                Root MSE          =     2500.8

------------------------------------------------------------------------------
             |               Robust
   price_res |      Coef.   Std. Err.      t    P>|t|     [95% Conf. Interval]
-------------+----------------------------------------------------------------
  weight_res |   1.823366   .7754338     2.35   0.021      .277567    3.369165
       _cons |   2.70e-06   290.7151     0.00   1.000      -579.53      579.53
------------------------------------------------------------------------------

. display .7808755*sqrt(71/72)
.77543379
\end{verbatim}
}
\onehalfspacing

The following code verifies the conclusion about the LZ standard errors. 
\singlespacing
{\footnotesize 
\begin{verbatim}
. replace rep78 = 0 if missing(rep78) 
(5 real changes made)

. reg price weight displ, cluster(rep78)

Linear regression                               Number of obs     =         74
                                                F(2, 5)           =       4.56
                                                Prob > F          =     0.0746
                                                R-squared         =     0.2909
                                                Root MSE          =     2518.4

                                  (Std. Err. adjusted for 6 clusters in rep78)
------------------------------------------------------------------------------
             |               Robust
       price |      Coef.   Std. Err.      t    P>|t|     [95% Conf. Interval]
-------------+----------------------------------------------------------------
      weight |   1.823366    .900214     2.03   0.099    -.4907079     4.13744
displacement |   2.087054   9.027184     0.23   0.826    -21.11806    25.29217
       _cons |    247.907   2043.732     0.12   0.908    -5005.675    5501.489
------------------------------------------------------------------------------

. regress price_res weight_res, cluster(rep78)

Linear regression                               Number of obs     =         74
                                                F(1, 5)           =       4.16
                                                Prob > F          =     0.0969
                                                R-squared         =     0.0609
                                                Root MSE          =     2500.8

                                  (Std. Err. adjusted for 6 clusters in rep78)
------------------------------------------------------------------------------
             |               Robust
   price_res |      Coef.   Std. Err.      t    P>|t|     [95% Conf. Interval]
-------------+----------------------------------------------------------------
  weight_res |   1.823366   .8939407     2.04   0.097    -.4745817    4.121314
       _cons |   2.70e-06   279.1706     0.00   1.000    -717.6308    717.6308
------------------------------------------------------------------------------

. display .900214*sqrt(71/72)
.89394066
\end{verbatim}
}
\onehalfspacing

The following code verifies the conclusion about the Newey--West standard errors with three choices of lags. If lag equals zero, the standard errors reduce to the EHW standard errors. 
\singlespacing
{\footnotesize 
\begin{verbatim}
. generate t = _n

. quiet tsset t

. newey price weight displ, lag(0)

Regression with Newey-West standard errors      Number of obs     =         74
maximum lag: 0                                  F(  2,        71) =      14.44
                                                Prob > F          =     0.0000

------------------------------------------------------------------------------
             |             Newey-West
       price |      Coef.   Std. Err.      t    P>|t|     [95% Conf. Interval]
-------------+----------------------------------------------------------------
      weight |   1.823366   .7808755     2.34   0.022     .2663445    3.380387
displacement |   2.087054   7.436967     0.28   0.780    -12.74184    16.91595
       _cons |    247.907   1129.602     0.22   0.827    -2004.455    2500.269
------------------------------------------------------------------------------

. newey price_res weight_res, lag(0)

Regression with Newey-West standard errors      Number of obs     =         74
maximum lag: 0                                  F(  1,        72) =       5.53
                                                Prob > F          =     0.0214

------------------------------------------------------------------------------
             |             Newey-West
   price_res |      Coef.   Std. Err.      t    P>|t|     [95% Conf. Interval]
-------------+----------------------------------------------------------------
  weight_res |   1.823366   .7754338     2.35   0.021      .277567    3.369165
       _cons |   2.70e-06   290.7151     0.00   1.000      -579.53      579.53
------------------------------------------------------------------------------

. display .7808755*sqrt(71/72)
.77543379

. 
. newey price weight displ, lag(1)

Regression with Newey-West standard errors      Number of obs     =         74
maximum lag: 1                                  F(  2,        71) =      10.62
                                                Prob > F          =     0.0001

------------------------------------------------------------------------------
             |             Newey-West
       price |      Coef.   Std. Err.      t    P>|t|     [95% Conf. Interval]
-------------+----------------------------------------------------------------
      weight |   1.823366   .7726505     2.36   0.021     .2827446    3.363987
displacement |   2.087054   7.989353     0.26   0.795    -13.84326    18.01737
       _cons |    247.907   1174.841     0.21   0.833    -2094.659    2590.473
------------------------------------------------------------------------------

. newey price_res weight_res, lag(1)

Regression with Newey-West standard errors      Number of obs     =         74
maximum lag: 1                                  F(  1,        72) =       5.65
                                                Prob > F          =     0.0201

------------------------------------------------------------------------------
             |             Newey-West
   price_res |      Coef.   Std. Err.      t    P>|t|     [95% Conf. Interval]
-------------+----------------------------------------------------------------
  weight_res |   1.823366   .7672661     2.38   0.020     .2938489    3.352883
       _cons |   2.70e-06   344.7862     0.00   1.000    -687.3187    687.3187
------------------------------------------------------------------------------

. display .7726505*sqrt(71/72)
.76726611

. 
. newey price weight displ, lag(2)

Regression with Newey-West standard errors      Number of obs     =         74
maximum lag: 2                                  F(  2,        71) =       9.77
                                                Prob > F          =     0.0002

------------------------------------------------------------------------------
             |             Newey-West
       price |      Coef.   Std. Err.      t    P>|t|     [95% Conf. Interval]
-------------+----------------------------------------------------------------
      weight |   1.823366   .7414398     2.46   0.016     .3449771    3.301755
displacement |   2.087054   8.096786     0.26   0.797    -14.05748    18.23159
       _cons |    247.907    1167.36     0.21   0.832    -2079.742    2575.556
------------------------------------------------------------------------------

. newey price_res weight_res, lag(2)

Regression with Newey-West standard errors      Number of obs     =         74
maximum lag: 2                                  F(  1,        72) =       6.13
                                                Prob > F          =     0.0156

------------------------------------------------------------------------------
             |             Newey-West
   price_res |      Coef.   Std. Err.      t    P>|t|     [95% Conf. Interval]
-------------+----------------------------------------------------------------
  weight_res |   1.823366   .7362729     2.48   0.016     .3556328    3.291099
       _cons |   2.70e-06     375.75     0.00   1.000    -749.0439    749.0439
------------------------------------------------------------------------------

. display .7414398*sqrt(71/72)
.73627291
\end{verbatim}
}
\onehalfspacing

\end{document}